\documentclass[12pt,reqno]{amsart}

\usepackage{amsthm,amsmath,amssymb,amsfonts,amscd,verbatim,latexsym}
\newtheorem{Theorem}{Theorem}[section]

\newtheorem{Proposition}[Theorem]{Proposition}

\newtheorem{Remark}[Theorem]{Remark}

\newtheorem{Example}[Theorem]{Example}

\newcommand{\complex}{\mathbf C}
\renewcommand{\P}{\mathbb P}
\renewcommand{\O}{\mathcal O}
\newcommand{\C}{\mathcal C}
\newcommand{\E}{\mathcal E}
\newcommand{\N}{\mathcal N}
\newcommand{\T}{\mathbb T}
\newcommand{\wE}{\widetilde E}

\newcommand{\ux}{\mathbf x}
\newcommand{\uy}{\mathbf y}
\newcommand{\ra}{\rightarrow}

\newcommand{\lra}{\longrightarrow}

\newcommand{\lp}{(}
\newcommand{\rp}{)}

\newcommand{\demo}{\noindent {\sc Proof.}\;}

\begin{document}
\title{On the Wronskian combinants of binary forms} 
\author[Abdesselam and Chipalkatti]
{Abdelmalek Abdesselam and Jaydeep Chipalkatti}
\maketitle

\parbox{12cm}{ \small 
{\sc Abstract.} 
For generic binary forms $A_1,\dots,A_r$ of order $d$ we construct 
a class of combinants $C = \{\C_q: 0 \le q \le r, q \neq 1\}$, to be 
called the Wronskian combinants of the $A_i$. 
We show that the collection $C$ gives a projective imbedding of the 
Grassmannian $G(r,S_d)$, and as a corollary, 
any other combinant admits a formula as an iterated transvectant 
in the $C$. Our second main result characterizes those collections of 
binary forms which can arise as Wronskian combinants.
These collections are the ones such that an associated algebraic
differential equation has the maximal number of
linearly independent
polynomial solutions.
Along the 
way we deduce some identities which connect Wronskians with 
transvectants.}

\bigskip \bigskip 

\parbox{12cm}{\small
Mathematics Subject Classification (2000):\, 13A50. \\ 
Keywords: binary forms, transvectant, combinant, Wronskian,
Grassmannian, Pl{\"u}cker imbedding, algebraic differential equations.} 

\bigskip 

\section{Introduction}
This article extends some of the investigations
in \cite{CBez} to the case of several binary forms. 
We begin by recalling the classical notion of a {\sl combinant} 
of binary forms (see \cite[\S 250]{GrYo}). 
A summary of our results will appear in 
\S \ref{results} below after the required notation is available. 
\subsection{} \label{intro1} 
Let $A_1,\dots,A_r$ denote generic 
forms of order $d$ in the variables $\ux = \{x_1,x_2\}$ 
(assume $r \le d$). Write
\begin{equation} A_i = \sum\limits_{j=0}^d \, \binom{d}{j} \, 
a_{ij} \, x_1^{d-j} x_2^j, \qquad 
(1 \le i \le r), \label{Ai.expression} \end{equation}
where the $a_{ij}$ are independent indeterminates. Given a matrix 
$g = \left( \begin{array}{cc} \alpha & \beta \\ \gamma & \delta 
\end{array} \right)$ such that $\alpha \, \delta - \beta \, \gamma =1$, 
make substitutions 
\[ x_1 = \alpha \, x_1' + \beta x_2', \quad 
   x_2 = \gamma \, x_1' + \delta x_2'; 
\] 
and now define $a_{ij}'$ by forcing the equalities 
\[ \sum\limits_{j=0}^d \, \binom{d}{j} \, a_{ij} \, x_1^{d-j} \, x_2^j = 
   \sum\limits_{j=0}^d \, \binom{d}{j} \, 
a_{ij}'\, {x_1'}^{d-j} \, {x_2'}^j. \] 
A polynomial function 
$Q(\{a_{ij}\};x_1,x_2)$ is called a {\sl joint covariant} of the 
the $\{A_i\}$ if 
\[ Q(\{a_{ij}\};x_1,x_2) = 
Q(\{a_{ij}'\};x_1',x_2'), \] 
for every $g$. 

It is called a {\sl combinant} if the following additional 
condition is satisfied: given a matrix $M \in SL_r$, define 
new constants $b_{ij}$ via the matrix equality 
${\mathcal B} = M \, {\mathcal A}$, 
where ${\mathcal B} = [b_{ij}], {\mathcal A} = [a_{ij}]$. Then 
we should have an equality 
\[ Q(\{a_{ij}\};x_1,x_2) = 
Q(\{b_{ij}\};x_1,x_2), \] 
for every $M$. We say that $Q$ is of degree $m$ and order $n$, if 
it has total degree $m$ in the coefficients of each $A_i$ and 
total degree $n$ in $\ux$. By the first fundamental theorem, 
the coefficients of $Q$ can be written 
as degree $m$ forms in the $r \times r$ minors of the matrix 
$[a_{ij}]$. 

For instance, for $r=2$ the resultant $R(A_1,A_2)$ is a combinant 
of degree $d$ and order zero. The Jacobian 
$\left| \begin{array}{cc} 
\frac{\partial A_1}{\partial x_1} & \frac{\partial A_1}{\partial x_2} \\ 
\frac{\partial A_2}{\partial x_1} & \frac{\partial A_2}{\partial x_2} 
\end{array} \right|$ is a combinant of degree one and order 
$2 d-2$. 

For fixed $(r,d)$, the combinants define a ring $R$ bigraded by $m$ and $n$. 
The structure of this ring can be very involved, 
and it is concretely known only for a few small values 
of $r$ and $d$ (see e.g.~\cite{Gordan1,Meulien1}). 
Our objective, roughly speaking, is to construct distinguished elements 
of this ring $C_0,C_2,\dots,C_r$ (sic) 
which generate it (in a slightly extended sense, to be made precise later). 

\subsection{Preliminaries}
Throughout, the base field will be $\complex$. 
We will write $S_d$ for the space of order $d$ forms in $\ux$, which 
is naturally a representation of $SL_2$. See~\cite[Ch.~11]{FH} 
or~\cite[Ch.~4]{Sturmfels} for standard facts about 
$SL_2$-representations. 
All of our constructions (and morphisms) will be $SL_2$-equivariant. 
Each finite-dimensional $SL_2$-module is 
canonically self-dual, we will use this identification 
without further comment. 

Our basic reference for invariant theory is~\cite{GrYo}; much of 
the same material is covered in~\cite{Glenn}. We will 
freely use the classical symbolic calculus interpreted {\it according to}
\cite[\S 2]{AC2} . 
Other treatments of this calculus,
which we will {\it not} use, can be
found in~\cite{KungRota,Olver}. 

If $E,F$ are two binary forms of orders $e,f$, their $k$-th transvectant 
is defined as 
\begin{equation}
(E,F)_k = \frac{(e-k)!(f-k)!}{e!f!}  \, 
\sum\limits_{i=0}^k \, (-1)^i \binom{k}{i} \, 
\frac{\partial^{\,k} E}{\partial x_1^{k-i} \, \partial x_2^i } \; 
\frac{\partial^{\,k} F}{\partial x_1^i \, \partial x_2^{k-i}} 
\; . \label{trans.defn} \end{equation}
It is identically zero outside the range $0 \le k \le \min \{e,f\}$. 

\subsection{} \label{comb1}
Let $G = G(r,S_d)$ denote the Grassmann variety 
of $r$-dimensional subspaces of $S_d$. 
Let $\lambda_m$ denote the partition 
$(\underbrace{m,\dots,m}_{\text{$r$ times}})$, and 
$S_{\lambda_m}$ the associated Schur functor (see \cite[Ch.~6]{FH}). 
By the Borel-Weil-Bott theorem (see~\cite[p.~687]{Porras})
we have an isomorphism of $SL_2$-representations 
\[ H^0(G,\O_G(m)) \simeq S_{\lambda_m}(S_d). \] 
An element of $S_{\lambda_m}(S_d)$ can be seen as a degree $m$ function 
in the Pl{\"u}cker coordinates via the imbedding 
\[ S_{\lambda_m}(S_d) \subseteq S_m(\wedge^r \, S_d). \] 
Then, a combinant $Q$ of degree $m$ and order $n$ 
can be identified (up to a scalar) with a morphism 
\[ \mu: S_n \lra H^0(G,\O_G(m)) \] 
which sends $F \in S_n$ to the transvectant $(F,Q)_n$. 
In the reverse direction, 
$\mu$  gives rise to a morphism 
\[ \mu': \complex \lra H^0(G,\O_G(m)) \otimes S_n, \] 
and then $Q$ can be recovered (up to scalar) as the element 
$\mu'(1)$. 

Hence combinants of degree-order $(m,n)$ (up to scalars) are in 
bijection with nonzero (and hence necessarily injective) morphisms 
$S_n \lra S_{\lambda_m}(S_d)$. 
In particular all the linear (i.e., degree one) combinants 
correspond to the irreducible summands of 
\[ \wedge^r S_d \simeq S_r(S_{d-r+1}). \] 
\begin{Example} \rm 
For $r=2$ and $m=1$, we have 
\[ S_2(S_{d-1}) \simeq
\bigoplus\limits_{i=1}^{\lfloor\frac{d+1}{2}\rfloor} \, 
S_{2d-2(2i-1)}. \] 
The $i$-th summand corresponds to the combinant $(A_1,A_2)_{2i-1}$. 
\end{Example} 
\begin{Example} \rm 
Assume $(r,d)=(2,5),m=2$. We have a plethysm decomposition 
\[ S_{\lambda_2}(S_5) = S_{(2,2)}(S_5) = 
S_{16} \oplus S_{12}^2 \oplus S_{10} \oplus 
S_8^3 \oplus S_6 \oplus S_4^3 \oplus S_0^2. \] 
(This was calculated using the Maple package `SF'.) 
This implies for instance, that two binary quintics $A_1,A_2$ 
have a two dimensional space of combinants of degree $2$ and order $12$. 
Writing $t_i = (A_1,A_2)_i$, a basis for this space is 
given by $(t_1,t_1)_2$ and $t_1 \, t_3$. 
\end{Example} 
Since the algebra $\bigoplus\limits_m S_{\lambda_m}(S_d)$ is generated in 
degree one, every combinant can be written as an iterated 
transvectant expression in linear combinants. 

\subsection{Wronskians} 
Given binary forms $F_1,\dots,F_s$ of order $n$, we define their Wronskian 
\[ W(F_1,\dots,F_s) = 
\left(\frac{(n-s+1)!}{n!}\right)^s \times \det 
\left(\frac{\partial^{s-1} F_i}{\partial x_1^{s-j} \, \partial x_2^{j-1}} 
\right)_{1 \le i,j \le s.} \] 
It is zero iff the $F_i$ are linearly dependent over $\complex$.
Using the classical symbolic calculus {\it according to \cite[\S 2]{AC2}},
if $F_i = {f^{(i)}_\ux}^n$ then 
\[ W = \prod\limits_{1 \le i<j \le s} 
(f^{(i)} \, f^{(j)}) \prod\limits_{1 \le i \le s} \, 
{f^{(i)}_\ux}^{n-s+1} \; . \] 
(The proof is easy: differentiate the symbolic expressions, 
and calculate the Vandermonde determinant.)
\subsection{Polarization} 
Introduce new letters $\uy = (y_1,y_2)$. If $\E$ is a form of 
order $n$ in $\ux$, then define its $k$-th polarization 
\[ \E^{\langle k \rangle} = \frac{(n-k)!}{n!} \, 
(y_1 \frac{\partial}{\partial x_1} + y_2 \frac{\partial}{\partial x_2})^k \, \E. \] 
\subsection{A summary of results} \label{results}
In this paper we will construct a set of linear combinants 
$C = \{C_0,C_2, \dots, C_r\}$ 
associated with a set of $r$ binary $d$-ics
$A_1,\dots,A_r$. In fact $C_0$ is the Wronskian of 
the $\{A_i\}$; the others are defined as transvectants of 
certain symbolic products derived from $\{A_i\}$. 
By construction $C_q$ is of order $r(d-r+1)-2q$. 
We will show that the $C$ enjoy a cluster of special properties: 
\begin{itemize} 
\item 
The subspace spanned by the $\{A_i\}$ can be recovered from 
$C$ as the solution space of the differential equation 
$\sum\limits_q \, (C_q,F)_{r-q} = 0$. 
\item 
The assignment 
$\text{span} \, \{A_i\} \lra [C_0,\dots,C_r]$ gives a projective 
imbedding of the Grassmannian $G(r,S_d)$. 
\item 
Every combinant $Q$ admits a formula of the type 
\[ Q = 
\frac{1}{C_0^N} \times 
(\text{A compound transvectant expression in the $C$}),
\]  
for some nonnegative integer $N$. 
\end{itemize} 

In Theorem~\ref{main2} we characterize
all possible values of $C$. Specifically we prove that a 
sequence of binary forms $E_0, E_2,\dots, E_r$ (of the correct orders) 
can arise as the Wronskian combinants of an $r$-dimensional subspace 
iff the differential equation 
\begin{equation} \sum\limits_q \, (E_q,F)_{r-q} = 0 
\label{diffeq.E} \end{equation}
admits $r$ linearly independent polynomial solutions. It follows that the 
image of the imbedding above is a determinantal variety, defined 
by equations of degree $d-r+2$. The proof proceeds in two steps: 
firstly we establish some identities which connect 
Wronskians with transvectants. Then we use these identities 
to `peel off' one summand at a time from equation~(\ref{diffeq.E}). 

\section{The Wronskian combinants} 
\label{construction.section}
Let us write $\N_r = \{q: 0 \le q \le r, q \neq 1\}$. 
In this section we will construct the Wronskian combinants 
$\C = \{C_q: q \in \N_r\}$. It is a special
instance of the method used by Gordan in order
to show that any covariant of a binary form is a linear combination
of iterated transvectants (see \cite[\S 2]{GordanBeweis}).
\subsection{} Let 
\[ A_i = {\alpha^{(i)}_{\ux}}^d, \quad 1 \le i \le r,\] 
denote $r$ binary $d$-ics, and $F=f_{\ux}^d$ another binary $d$-ic. Then 
\[ W = W(A_1,\dots,A_r,F) = f_{\ux}^{d-r} 
\prod\limits_{i=1}^r {\alpha_{\ux}^{(i)}}^{d-r} 
\prod\limits_{1 \le i < j \le r} (\alpha^{(i)} \, \alpha^{(j)}) 
\prod\limits_{i=1}^r (\alpha^{(i)} \, f).
\] 
We will rewrite this expression as follows. Define 
\[ \Phi = \prod\limits_{i=1}^r {\alpha_{\ux}^{(i)}}^{d-r} \, 
\prod\limits_{1 \le i < j \le r} (\alpha^{(i)} \, \alpha^{(j)}) \, 
\prod\limits_{i=1}^r \alpha^{(i)}_{\uy}. \] 
Then 
\begin{equation} 
W = \left. (\Phi,f_{\uy}^d)^{\uy}_r \right|_{\uy:=\ux}, 
\label{wphi} \end{equation} 
which is to say, take the $r$-th transvectant of the pair 
$\Phi,f_{\uy}^d$ as forms in $\uy$ (treating the $\ux$ as constants), 
and then substitute $\ux$ for $\uy$. 

If we write $G_{\ux} = \prod\limits_{i=1}^r \alpha_{\ux}^{(i)}$, then 
\[ \Phi = \{ \, \prod\limits_{1\le i<j\le r} 
(\alpha^{(i)} \, \alpha^{(j)})\} \, G_\ux^{d-r} G_\uy. 
\]
Now the Gordan series (see \cite[\S 52]{GrYo}) gives an identity 
\[ G_\ux^{d-r} G_\uy = \sum\limits_{q=0}^r \, 
\frac{\binom{rd-r^2}{q} \binom{r}{q}}{\binom{rd-r^2+r-q+1}{q}} \, 
[(G_\ux^{d-r},G_\ux)_q]^{\langle r-q \rangle} \, (\ux \, \uy)^q. \] 
Hence let us define 
\begin{equation} 
C_q = (-1)^q \, 
\frac{\binom{rd-r^2}{q} \binom{r}{q}}{\binom{rd-r^2+r-q+1}{q}} \; 
\{ \prod\limits_{1 \le i < j \le r} (\alpha^{(i)} \, \alpha^{(j)}) \} \; 
(G_\ux^{d-r},G_\ux)_q. \label{definition.Cq} \end{equation} 
By construction this is a combinant which is linear in the 
coefficients of each $A_i$, and of order $r(d-r+1)-2q$. 
Notice that $C_1$ is identically zero, because it is the Jacobian 
of two functionally dependent forms (and hence it will not be mentioned any 
further). Since $G_\ux$ is of order $r$, the index $q$ ranges over $\N_r$. 

\subsection{} Thus we have 
\[ \Phi = \sum\limits_{q \in \N_r} \, 
(-1)^q \, [C_q]^{\langle r-q \rangle} \, (\ux \, \uy)^q. \] 
We will calculate $W$ by substituting this expression into 
(\ref{wphi}). Write 
$C_q = c_\ux^{\, r(d-r+1)-2q}$. 
Then 
\[ [C_q]^{\langle r-q \rangle} \, (\ux \, \uy)^q = 
c_\ux^{\, r(d-r)-q} \, c_\uy^{r-q} \, (\ux \, \uy)^q, \]
and 
\[ ([C_q]^{\langle r-q \rangle} \, (\ux \, \uy)^q,f_\uy^d)_r^{\uy} = 
(-1)^q \, c_\ux^{r(d-r)-q} \,(c \, f)^{r-q} \, f_\ux^q \, f_\uy^{d-r}. \] 
Letting $\uy:=\ux$, this reduces to 
\[ (-1)^q \, c_\ux^{\, r(d-r)-q} \,(c \, f)^{r-q} \, f_\ux^{d-r+q} = 
(-1)^q \, (C_q,F)_{r-q}. \] 
In sum, 
\begin{equation} W = \sum\limits_{q \in \N_r} \, (C_q,F)_{r-q}. 
\label{eqn.W} \end{equation}
\subsection{} 
Given an arbitrary collection of forms $E = \{E_q: q \in \N_r\}$ of 
orders $r(d-r+1)-2q$, we define 
\[ \psi_E(F) = \sum\limits_{q \in \N_r} \, (E_q,F)_{r-q}. \] 
Then $\psi_E(F)=0$ is an 
algebraic differential equation of order $r$ dependent on 
the parameters $E$. 

Let $C = \{C_q\}$ be the combinants constructed above 
associated with $\{A_1,\dots,A_r\}$. If the $\{A_i\}$ are linearly 
dependent, then (and only then) all $\C_q$ are zero. 
If they are independent, then a binary $d$-ic $F$ belongs to their 
linear span iff $W(A_1,\dots,A_r,F)=0$, i.e., iff $\psi_C(F)=0$. 
Hence the $C$ completely characterize the subspace spanned by the $A$. 

\begin{Example} \rm 
Assume $r=2$, then $C_0 = (A_1,A_2)_1$. 
Since $G_\ux$ is a quadratic, 
\[ (G_\ux^{d-2},G_\ux)_2 = 
\frac{2-d}{2(2d-5)} \, (G_\ux, G_\ux)_2 \; G_\ux^{d-3}. \] 
(This can be checked by a direct symbolic calculation.) Hence 
\[ C_2 = \frac{\binom{2d-4}{2}}{\binom{2d-3}{2}} \, 
\frac{2-d}{2(2d-5)} \, 
(\alpha^{(1)} \, \alpha^{(2)})^3 \, 
{\alpha^{(1)}_\ux}^{d-3} \, {\alpha^{(2)}_\ux}^{d-3}  = 
\frac{2-d}{4d-6} \, (A_1,A_2)_3. \] 
\end{Example}

\section{The incomplete Pl{\"u}cker imbedding} \label{plucker}
Let 
\[ U = \bigoplus\limits_{q \in \N_r} \, S_{r(d-r+1)-2q} \, , \] 
and consider the morphism 
\[ \pi: G(r,S_d) \lra \P \, U\] 
which sends the subspace $\Lambda = \text{span}\{A_1,\dots,A_r\}$ to 
$[C_0,C_2,\dots,C_r]$. 
\begin{Theorem} \sl 
The morphism $\pi$ is an imbedding. 
\end{Theorem} 
\demo 
Since $\Lambda$ can be recovered from $C$, we deduce that 
$\pi$ is set-theoretically injective. It remains to show that 
$\pi$ is injective on tangent spaces (\cite[Ch.~II, Prop.~7.3]{Ha}). 

The tangent space $T_{G,\Lambda}$ is canonically 
isomorphic to $\text{Hom}(\Lambda,S_d/\Lambda)$. 
Assume that $v \in T_{G,\Lambda}$ sends 
$A_i$ to $B_i + \Lambda$. 

The tangent space to $\P \, U$ at $\pi(v)$ is isomorphic to 
$U/[C_0,\dots,C_q]$. To calculate the image vector $d \pi(v)$, 
define \[ D_q = 
\lim\limits_{\epsilon \ra 0} \, 
\frac {C_q(A_1 + \epsilon \, B_1,\dots,A_r + \epsilon \, B_r)-
C_q(A_1,\dots,A_r)}{\epsilon}. \] 
Then $d \pi(v) = [D_0,D_2,\dots,D_q]$, considered modulo 
$[C_0,\dots,C_q]$. 
The Wronskian combinants are multilinear in each argument, hence 
\[ D_q = 
\sum\limits_{i=1}^r \, C_q(A_1,\dots,{\widehat A_i}|B_i, \dots, A_r),
\] 
where the last expression means that $A_i$ is to be replaced by $B_i$. 
Assume that $d \pi(v)=0$, so 
there exists a constant $\alpha$ such that $D_q = \alpha \, C_q$ 
for all $q$. But then $\sum\limits_q \, (D_q,A_1)_{r-q}=0$, i.e., 
\[ \sum\limits_{i=1}^r \; 
\left\{ \sum\limits_q \; 
(C_q(A_1,\dots,{\widehat A_i}|B_i, \dots, A_r),A_1)_{r-q} \right\} 
=0. \] 
All the summands except $i=1$ vanish for obvious reasons, hence 
so does the remaining one. This implies that 
$A_1 \in \text{Span} \, \{B_1,A_2,\dots,A_r\}$, 
which forces $B_1 \in \text{Span} \{A_1,\dots,A_r\}$. 
Similarly each $B_i \in \Lambda$, hence $v$ must be the 
zero vector. This shows that $\pi$ was injective on tangent 
spaces. The theorem is proved. \qed 

\subsection{} 
Let $Q$ be an arbitrary combinant of $r$ binary $d$-ics. 
We will show that $Q$ admits a `formula' as mentioned in the 
introduction. In order to make this precise, assume 
$A_i$ to be as in equation~(\ref{Ai.expression}). 
Let $\T$ denote the smallest $\complex$-subalgebra of the polynomial algebra 
$\complex [\{a_{ij}\}_{i,j},\ x_1,x_2]$ such that 
\begin{itemize} 
\item 
$C_0,\dots,C_q \in \T$, and 
\item 
if $e_1,e_2 \in \T$, then $(e_1,e_2)_k \in \T$ for all $k \ge 0$. 
\end{itemize}
Each element of $\T$ is a combinant, and there is a natural bigraded 
decomposition of $\T$ induced by the degree $m$ and order $n$. 
For instance, the element 
\[ ((C_0,C_3)_3, C_2)_2 + 5 \, (C_0^2,C_4)_6 \] 
is bihomogeneous of degree $3$ and order $3 \, r(d-r+1)-20$. 

\begin{Theorem} \sl Let $Q$ be a combinant of the $\{A_i\}$. 
Then there exists an integer $N \ge 0$ such that 
$C_0^N \, Q \in \T$. 
\label{formula.Q} \end{Theorem} 
Since $C_0$ is always nonzero on linearly independent forms, this 
shows the existence of a formula for $Q$. 

\smallskip 

\demo 
Given the imbedding $G \subseteq \P \, U$, for every integer $m \ge 0$ 
we have the restriction morphism 
\[ f_m: 
\underbrace{H^0(\P \, U,\O_{\P}(m))}_{= \, S_m(U)}
\lra H^0(G,\O_G(m)). \] 
Let $m = N + \deg Q$. 
The combinant $C_0^N \, Q$ will lie in $\T$ iff the image 
of the corresponding morphism (see~\S\ref{comb1}) 
\[ S_{\text{ord}(C_0^N Q)} \lra H^0(G,\O_G(m)) \] 
is contained in the image of $f_m$. 
But this can always be arranged by choosing $N > > 0$, 
since $f_m$ is surjective for $m > > 0$. \qed 

\begin{Remark} \rm 
If $\rho$ is the Castelnuovo 
regularity of the ideal sheaf ${\mathcal I}_G$, then 
$H^1(\P U,{\mathcal I}_G(\rho-1))=0$ implying that 
$f_{\rho-1}$ is surjective. Hence $N$ can be chosen to be 
$\max \, \{0,\rho - \deg Q-1\}$. It is possible (but rather tedious) 
to calculate an explicit upper bound for $\rho$ 
from the Hilbert polynomial of 
$G$ (see~\cite{KL.SGA6}), but we have not attempted this. 
\end{Remark} 

\begin{Example} \rm 
Assume $(r,d)=(2,5)$, and write $t_i = (A_1,A_2)_i$. Then 
the element $C_0 \, t_5$ lies in $\T$, in fact 
there is an identity 
\[ t_5 = \frac{1}{C_0} \, 
[ 50 \, C_2^2 -15 \, (C_0,C_0)_4 -40 \, (C_0,C_2)_2 ]. \] 
{\sc Proof sketch:} $C_0 \, t_5$ is of degree-order 
$(2,8)$. The plethysm 
$S_{(2,2)}(S_5)$ contains $3$ copies of $S_8$, hence there is 
a $3$ dimensional space of such combinants. By specialising $A_1,A_2$
we can show that $C_2^2,(C_0,C_0)_4,(C_0,C_2)_2$ 
are linearly independent, hence they form a basis of this space. Thus 
$C_0 \, t_5$ must be 
expressible as their linear combination. To find the actual coefficients 
we only need to solve a system of linear equations. 
\end{Example}

\section{Wronskians and Transvectants} 
We now come to our second main theorem which characterizes all 
possible values of $C$. 
\subsection{} Let $E = \{E_q: q \in \N_r \}$ be an arbitrary 
collection of binary forms of orders $r(d-r+1)-2q$, such 
that $E_1 \neq 0$. 
In general the $r$-th order differential equation 
\begin{equation} \psi_E(F)=0 \label{deq.E} \end{equation}
does not admit $r$ linearly independent polynomial solutions. Our 
result says that if indeed it does, then the $E$ must be 
values of Wronskian combinants. 
\begin{Theorem} \sl 
Assume that there exist $r$ linearly independent 
$d$-ics $A_1,\dots,A_r$ such that 
$\psi_E(A_i)=0$. Then there exists a nonzero constant $k$ such that 
\[ E_q = k \, C_q(A_1,\dots,A_r) \] 
for all $q \in \N_r$. 
\label{main2} \end{Theorem} 
Then, of course, we can arrange that $E_q = C_q$ by replacing $A_1$ with 
$k \, A_1$. 
\subsection{} 
The proof hinges upon certain identities 
involving transvectants and Wronskians. 
Let $B$ denote a form of order $n$. 
For $0 \le p \le \min \, \{d,n\}$, define 
\[
\Gamma_p(B; A_1,\dots,A_r) = \sum\limits_{i=1}^r
(-1)^{i+1} \, (B,A_i)_p \, 
W(A_1,\ldots, {\widehat{A_i}},\ldots,A_r). 
\]
We will tentatively abbreviate this to $\Gamma_p$. 
Now the key result is the following.
\begin{Proposition} \sl
We have identities 
\[ \Gamma_p = 
\begin{cases} 
\qquad 0 & \text{for $0 \le p \le r-2$}, \\ 
(-1)^{r-1} \, B \, W(A_1,\ldots,A_r) & 
\text{for $p = r-1$}, \\ 
(-1)^{r-1} \, r \, (B,W(A_1,\ldots,A_r))_1 
& \text{for $p = r$.}
\end{cases}\] 
\label{keyprop} \end{Proposition} 
We have found no such simple identities for $p > r$. The proposition will be 
proved in~\S \ref{proof.keyprop}; 
meanwhile let us use it to prove the theorem. 

\medskip 

\noindent {\sc Proof of theorem~\ref{main2}.} 
By hypothesis
\begin{equation} 
\sum\limits_{q \in \N_r} \, (E_q, A_i)_{r-q}=0.  \label{psiE}
\end{equation} 
Multiply this equation by 
$(-1)^{i+1} \, W(A_1,\ldots, {\widehat{A_i}},\ldots,A_r)$ 
and sum over $1 \le i \le r$. This gives 
\[ \sum\limits_{q \in \N_r}
\Gamma_{r-q}(E_q; A_1,\ldots, A_r) =0. \] 
By the proposition, we have 
$\Gamma_{r-q}=0$ for $r-q\le r-2$, i.e., for  $q \ge 2$.
\emph{Fortunately} there is no $q=1$ term, hence 
\[ \Gamma_r(E_0; A_1,\ldots, A_r)= 
(E_0,W(A_1,\ldots, A_r))_1 =0. \] 
In general, if $M,N$ are forms of the same order, then $(M,N)_1=W(M,N)$; 
which can be zero only if $M,N$ are multiples of each other. Hence 
there exists a constant $k$ such that
\[ E_0= k \, W(A_1,\ldots, A_r) = k \, C_0. \] 
Now write $\wE_q = E_q - k \, C_q$. Subtract the equation 
$k \, \psi_C(A_i)=0$ from~(\ref{psiE}), this gives 
\begin{equation} 
\sum\limits_{q=2}^r \, (\wE_q, A_i)_{r-q}=0. 
\label{descentready} \end{equation} 
Multiply~(\ref{descentready}) by 
$(-1)^{i+1} \, W(A_1,\ldots, {\widehat{A_{i}}},\ldots,A_{r-1})$, 
(note that $A_r$ is missing), 
and sum over $1 \le i \le r-1$. Then we have 
\[ \sum\limits_{q=2}^r
\Gamma_{r-q}(\wE_q; A_1,\ldots, A_{r-1})=0. 
\]
By the proposition, all the summands for 
$r - q \le r-3$, (i.e., $q \ge 3$) are zero. Hence 
\[
\Gamma_{r-2}(\wE_2; A_1,\ldots, A_{r-1})= 
(-1)^{r-2} \, 
\wE_2 \, W(A_1,\ldots,A_{r-1})=0. \] 
Since the $A_i$ are linearly independent,  
the Wronskian on the right is nonzero, hence $\wE_2=0$.
We can now repeat this procedure by dropping $A_{r-1},A_{r-2}$ etc.; 
this will successively force $\wE_3,\wE_4$~etc.~to be zero. \qed 

\subsection{} 
The assignment $F \lra \psi_E(F)$ gives a morphism 
\[ \psi: S_d \otimes \O_{\P U}(-1) \lra 
S_{d + r(d-r-1)} \otimes \O_{\P U}. \] 
The theorem implies that $\pi(G)$ is equal to the locus 
$\{ \text{rank} \, \psi \le d-r+1 \}$. 
Consequently $\pi(G)$ is set-theoretically 
defined by equations of degree $d-r+2$. 

\section{Proof of proposition~\ref{keyprop}.} 
\label{proof.keyprop}
\subsection{} 
Let us write $B=\beta_\ux^n$. Then $\Gamma_p$ has the symbolic expression
\[ \sum\limits_{i=1}^r \; (-1)^{i+1} \, 
\{(\beta \, \alpha^{(i)})^p \, 
\beta_\ux^{n-p}\,  {\alpha_{\ux}^{(i)}}^{d-p} 
\prod\limits_{{1\le j<k\le r}\atop{j,k\neq i}} 
(\alpha^{(j)}\, \alpha^{(k)}) 
\prod\limits_{{1\le j\le r}\atop{j\neq i}} 
{\alpha_{\ux}^{(j)}}^{d-r+2} \}. \]
Now dehomogenize using the following substitutions: 
\[ (\beta_1,\beta_2)=(b,1), \quad 
(\alpha_{1}^{(i)},\alpha_{2}^{(i)})=(a_i,1), \quad 
(x_1,x_2)=(1,-u).  \]
Then we have 
\[\begin{aligned} 
\Gamma_p = & \, \sum\limits_{i=1}^r \, (-1)^{i+1} \, \{
(b-a_i)^p \, (b-u)^{n-p} \, (a_i-u)^{d-p} \\ 
& \times 
\prod\limits_{{1\le j<k\le r}\atop{j,k\neq i}} (a_j-a_k) \times 
\prod\limits_{{1\le j\le r}\atop{j\neq i}} (a_j-u)^{d-r+2} \, \}.
\end{aligned} \] 
This can be rewritten as 
$(b-u)^{n-p} \times \prod\limits_{i=1}^r (a_i-u)^{d-r+2} \, \times$ 
\[
\sum\limits_{i=1}^r \, \left\{ (-1)^{i+1} \, 
(b-a_i)^p \, (a_i-u)^{r-p-2} \times
\left| \begin{array}{cccc}
a_1^{r-2} & \cdots & a_1 & 1 \\
\vdots & {\ } & \vdots & \vdots \\
a_{i-1}^{r-2} & \cdots & a_{i-1} & 1 \\
a_{i+1}^{r-2} & \cdots & a_{i+1} & 1 \\
\vdots & {\ } & \vdots & \vdots \\
a_r^{r-2} & \cdots & a_r & 1 
\end{array} \right| \; \right\}. \]
Hence 
\begin{equation} \Gamma_p=(b-u)^{n-p}\times
\prod\limits_{i=1}^r (a_i-u)^{d-r+2}
\times \left|
\begin{array}{ccccc}
Q(a_1) &  a_1^{r-2} & \cdots & a_1 & 1 \\
\vdots & \vdots & {\ } & \vdots & \vdots \\
Q(a_r) &  a_r^{r-2} & \cdots & a_r & 1 
\end{array} \right|,  \label{Gammap.Q} \end{equation}
with $Q(a)=(b-a)^p \, (a-u)^{r-p-2}$.
(To see this, expand the last determinant by its first column.) 

If $p\le r-2$, then $Q(a)$ is a polynomial in $a$ of 
degree $r-2$, hence the first column is a linear combination of the others.
This forces $\Gamma_p=0$, which is the first part of the proposition. 
\qed 

\subsection{} 
Now let $p=r-1$, so that
\[ \begin{aligned} 
Q(a) & =\frac{(b-a)^{r-1}}{a-u} = 
(-1)^{r-1} \, \frac{\left[ (a-u)+(u-b) \right]^{r-1}} {a-u} \\ 
& =(-1)^{r-1} \, \sum\limits_{j=0}^{r-1} 
\binom{r-1}{j} (a-u)^{j-1} (u-b)^{r-1-j}
\end{aligned} \]
By the previous argument on columns, we only need the $j=0$ term to 
calculate the determinant. Hence $\Gamma_{r-1}$ is equal to 
\[ \begin{aligned} 
{} & (b-u)^{n-r+1} 
(u-b)^{r-1} \times \prod\limits_{i=1}^r \, (a_i-u)^{d-r+2} \times 
\left|
\begin{array}{ccccc}
\frac{1}{a_1-u} &  a_1^{r-2} & \cdots & a_1 & 1 \\
\vdots & \vdots & {\ } & \vdots & \vdots \\
\frac{1}{a_r-u} &  a_r^{r-2} & \cdots & a_r & 1 
\end{array}
\right| \\ 
=  - & (b-u)^n \, \prod\limits_{i=1}^r
(a_i-u)^{d-r+2} \times 
\left| \begin{array}{ccccc}
\frac{1}{u-a_1} &  a_1^{r-2} & \cdots & a_1 & 1 \\
\vdots & \vdots & {\ } & \vdots & \vdots \\
\frac{1}{u-a_r} &  a_r^{r-2} & \cdots & a_r & 1 
\end{array} \right| \end{aligned} \] 

The last determinant can be written as 
\[ \sum\limits_{s \ge 0} \; 
\frac{1}{u^{s+1}} \left|
\begin{array}{ccccc}
a_1^s &  a_1^{r-2} & \cdots & a_1 & 1 \\
\vdots & \vdots & {\ } & \vdots & \vdots \\
a_r^s &  a_r^{r-2} & \cdots & a_r & 1 
\end{array} \right| \;. \] 

\medskip 

Let us write $\Delta_r(a)$ for $\prod\limits_{1\le j<k\le r}
(a_j-a_k)$. 
Using the Schur polynomials $S_\lambda(a)$ 
(see \cite[Appendix 1]{FH}), we rewrite $\Gamma_{r-1}$ as 
\[ \begin{aligned} 
{} & - (b-u)^n \, 
\{ \prod\limits_{i=1}^r \, 
(a_i-u)^{d-r+2} \} \times \Delta_r(a) \, 
\sum\limits_{s \ge r-1} \; \frac{1}{u^{s+1}} \, 
S_{(s-r+1)}(a_1,\dots,a_r) \\ 
= & - (b-u)^n \, 
\{ \prod\limits_{i=1}^r \, (a_i-u)^{d-r+2} \} \times
\Delta_r(a) \times
\frac{1}{u^r}\times
\prod\limits_{i=1}^r\frac{1}{1-\frac{a_i}{u}} \\ 
= & \, (-1)^{r-1} \, (b-u)^n \, \Delta_r(a) \, 
\prod\limits_{i=1}^r \, (a_i-u)^{d-r+1}. 
\end{aligned} \]
Rehomogenizing this, we get 
\[ \Gamma_{r-1}=
(-1)^{r-1} \, \beta_\ux^n
\times
\prod\limits_{i=1}^r \, {\alpha_{\ux}^{(i)}}^{d-r+1}
\times 
\prod\limits_{1 \le j < k \le r} \, (\alpha^{(j)} \, \alpha^{(k)}), 
\]
which proves the second part. 

\subsection{} 
Finally let $p=r$. Then $Q(a)$ is equal to 
\[ \begin{aligned} 
{} & =\frac{(b-a)^r}{(a-u)^2}= (-1)^r \, 
\frac{[\, (u-b)+(a-u) \,]^r}{(a-u)^2} \\ 
& =(-1)^r \left[
\frac{(u-b)^r}{(a-u)^2}+ \frac{r \, (u-b)^{r-1}}{a-u}+ 
\text{irrelevant terms} \right]
\end{aligned} \] 
Now substitute this into~(\ref{Gammap.Q}). The positive powers of 
$(a-u)$ contribute nothing to the sum, hence 
\[ \begin{aligned} 
\Gamma_r  =(-1)^r \, & (b-u)^{n-r} \, \{ \, \prod\limits_{i=1}^r \, 
(a_i-u)^{d-r+2} \, \} \\ 
\times \, & \left[(u-b)^r \, D_2 + r \, (u-b)^{r-1} \, D_1\right], 
\end{aligned} \] 
where
\[ D_\nu= \left|
\begin{array}{ccccc}
\frac{1}{(a_1-u)^\nu} &  a_1^{r-2} & \cdots & a_1 & 1 \\
\vdots & \vdots & {\ } & \vdots & \vdots \\
\frac{1}{(a_r-u)^\nu} &  a_r^{r-2} & \cdots & a_r & 1 
\end{array} \right|.\]
As in the previous case, 
\[ D_1=-\frac{\Delta_r(a)}{u^r} \times
\prod\limits_{i=1}^r \, \frac{1}{1-\frac{a_i}{u}}
=\frac{(-1)^{r-1}\Delta_r(a)}{\prod\limits_{i=1}^r (a_i-u)}, 
\] whereas
\[ D_2=\frac{\partial D_1}{\partial u}= D_1
\times \sum\limits_{i=1}^r \, \frac{1}{a_i-u}. \]

\smallskip 

Putting everything together, $\Gamma_r$ equals 
\[ \begin{aligned} 
{} & (b-u)^n \, \{ \prod\limits_{i=1}^r
(a_i-u)^{d-r+2} \} \times \left[D_2+\frac{r\ D_1}{u-b} \right] \\ 
= & \, (-1)^{r-1} (b-u)^n\ \Delta_r(a) \, 
\{ \prod\limits_{i=1}^r \, (a_i-u)^{d-r+1} \} 
\times \sum\limits_{i=1}^r 
\left\{ \frac{1}{a_i-u}+ \frac{r}{u-b} \right\} \\ 
= & \, (-1)^{r-1} \, (b-u)^{n-1} \, \Delta_r(a) \, 
\{ \, \prod\limits_{i=1}^r \, (a_i-u)^{d-r+1} \, \} 
\times \sum\limits_{i=1}^r \, 
\frac{b-a_i}{a_i-u}. \end{aligned} \] 
Homogenizing, 
\[ \begin{aligned} 
\Gamma_r & =(-1)^{r-1} \, \beta_\ux^{n-1} \, 
\times\prod\limits_{1\le j<k\le r} (\alpha^{(j)} \, \alpha^{(k)}) \\ 
& \times \sum\limits_{i=1}^r \, 
\left((\beta \, \alpha^{(i)})
\times {\alpha_\ux^{(i)}}^{d-r}\times
\prod\limits_{{1\le j\le r}\atop{j\neq i}}
{\alpha_\ux^{(j)}}^{d-r+1}
\right) \\ 
& = (-1)^{r-1} \, r \, {\lp B,W(A_1,\ldots,A_r) \rp}_1. 
\end{aligned} \] 
In the last step we have used the general formula for 
transvectants of symbolic products (see~\cite[\S 3.2.5]{Glenn}). 
The proposition is proved; 
this also completes the proof of Theorem~\ref{main2}. \qed 

\smallskip 

The following question arises naturally: given a reductive group $G$ and 
a $G$-module $V$, investigate how much of the theory carries over to 
the Grassmannian $G(r,V)$. 

\smallskip 

{\sc Acknowledgements.} {\small The first author is grateful to
D.~Brydges and J.~Feldman for their invitation to the
University of British Columbia.
The second author is grateful for the 
financial support by NSERC. The `SF' (Symmetric Functions) package for 
Maple has been useful, and we are grateful to its author John Stembridge. 
The following electronic libraries have been useful in accessing classical
references:

\begin{itemize}

\item The G\"ottinger DigitalisierungsZentrum ({\bf GDZ})

\item Project Gutenberg ({\bf PG})

\item The University of Michigan Historical Mathematics Collection ({\bf UM})

\end{itemize}

\bibliographystyle{plain}

\bigskip
\centerline{------------------------------}

\vspace{3cm}

\parbox{6cm}{\small 
{\sc Abdelmalek Abdesselam} \\ 
Department of Mathematics\\
University of British Columbia\\
1984 Mathematics Road \\
Vancouver, BC V6T 1Z2 \\ Canada. \\ 
{\tt abdessel@math.ubc.ca}} 
\hfill 
\parbox{5cm}{\small 
LAGA, Institut Galil\'ee \\ CNRS UMR 7539\\
Universit{\'e} Paris XIII\\
99 Avenue J.B. Cl{\'e}ment\\
F93430 Villetaneuse \\ France.} 

\vspace{1.5cm} 

\parbox{6cm}{\small 
{\sc Jaydeep Chipalkatti} \\ 
Department of Mathematics\\
University of Manitoba \\ 
433 Machray Hall \\ 
Winnipeg MB R3T 2N2 \\ Canada. \\ 
{\tt chipalka@cc.umanitoba.ca}}

\end{document}